\documentclass{amsart}
\usepackage{amsmath,amsthm,amsfonts,amssymb,xcolor}
\usepackage{mathrsfs}
\usepackage{enumitem}
\usepackage[colorlinks,linkcolor=blue,anchorcolor=blue,citecolor=blue,urlcolor=blue]{hyperref}
\usepackage[initials,nobysame]{amsrefs}

\renewcommand{\PrintDOI}[1]{DOI \href{https://doi.org/#1}{\nolinkurl{#1}}}

\makeatletter
\@namedef{subjclassname@2020}{\textup{2020} Mathematics Subject Classification}
\makeatother
\newtheorem{thm}{Theorem}[section]
\newtheorem{lem}[thm]{Lemma}

\newtheorem{cor}[thm]{Corollary}
\theoremstyle{definition}

\theoremstyle{remark}
\newtheorem{rem}[thm]{Remark}
\DeclareMathOperator{\Var}{Var}
\DeclareMathOperator{\osc}{osc}
\DeclareMathOperator{\diam}{diam}
\DeclareMathOperator{\sgn}{sgn}
\newcommand{\R}{\mathbb R}
\newcommand{\1}{\mathbf 1}
\newcommand{\C}{\mathcal C}

\allowdisplaybreaks[1]
\numberwithin{equation}{section}

\hypersetup{
  pdftitle={Variational Estimates for the Hardy--Littlewood Maximal Operator on Finite Graphs},
  pdfauthor={Shengzhe Lu and Fayou Zhao}
}

\begin{document}
\title[Variational estimates on finite graphs]{Variational Estimates for the Hardy--Littlewood Maximal Operator on Finite Graphs}
\author[S. Lu and F. Zhao]{Shengzhe Lu and Fayou Zhao$^*$}
\address[Shengzhe Lu]{Department of Mathematics, Shanghai University, Shanghai 200444, P.R. China}
\email{lushengzhe@shu.edu.cn}
\address[Fayou Zhao]{Newtouch Center for Mathematics of Shanghai University, Department of Mathematics, Shanghai University, Shanghai 200444, P.R. China}
\email{fyzhao@shu.edu.cn}
\thanks{$\ast$ Corresponding author}
\subjclass[2020]{Primary 42B25; Secondary 26A45, 42B35, 05C12.}
\keywords{Hardy--Littlewood maximal operator, finite graph, edge variation, sharp constant, diameter}
\date{}
\begin{abstract}
We study edge variation for the centered Hardy--Littlewood maximal operator on finite connected graphs. We determine the optimal uniform coefficient controlling the oscillation of the maximal function by the edge variation of the original function and derive two-sided estimates for the variational constant of an individual graph. These estimates involve point-mass tests, the diameter, and the number of edges. On star graphs, we determine the sharp constant for exponents greater than one and characterize all nonconstant real-valued extremizers for every exponent.
\end{abstract}
\maketitle

\section{Introduction}\label{sec:intro}
Kinnunen \cite{Kinnunen} proved that the Hardy--Littlewood maximal operator is bounded on $W^{1,p}(\R^d)$ for $1<p\leq\infty$. On a finite graph, differences across edges provide a discrete counterpart of derivatives. This leads to the problem of determining the smallest constant in the inequality
\[
\Var_p(M_Gf)\leq\C_{G,p}\Var_p(f),
\qquad 0<p\leq\infty,
\]
where $M_G$ is the centered Hardy--Littlewood maximal operator and $\Var_p$ denotes edge variation. The definitions used here count each undirected edge once and take averages with respect to counting measure.

Liu and Xue \cite{LX} studied variational inequalities on finite graphs. Gonz\'alez-Riquelme and Madrid \cite{GRM1} established further estimates on complete graphs, star graphs, and general finite graphs. In particular, they determined the sharp constants on complete and star graphs for several ranges of exponents. Their subsequent paper \cite{GRM2} settled the case $p=2$ on star graphs. More recently, Gonz\'alez-Riquelme, Kova\v{c}, and Madrid \cite{GKM} obtained the sharp constant on complete graphs for every finite positive exponent and proved that $\C_{G,p}$ is unbounded over all finite connected graphs for each fixed $0<p<\infty$.

Our first result gives the optimal uniform estimate of $\osc(M_Gf)$ in terms of $\Var_p(f)$ over all connected graphs with a prescribed number of vertices. More precisely, Theorem~\ref{thm:sharp-osc-variation} determines $\Omega_p(n)$ such that
\[
\sup_{G\in\mathscr G_n}\ \sup_{\Var_p(f)>0}
\frac{\osc(M_Gf)}{\Var_p(f)}=\Omega_p(n),
\]
where $\mathscr G_n$ is the class of connected simple undirected graphs on $n$ vertices. Equality is attained on a path for every exponent. This oscillation estimate, combined with point-mass tests and a bound for edge variation in terms of oscillation, yields two-sided estimates for $\C_{G,p}$ that retain the actual diameter. The edge-count estimate used in the argument is due to Ore \cite{Ore}. Maximization over the possible diameters then gives bounds depending only on $n$ and $p$.

Theorem~\ref{thm:star} determines $\C_{S_n,p}$ for every $0<p\leq\infty$. For $0<p\leq1$, the sharp constant is $(n-1)/n$. For $p>1$, the formula has two branches, separated, when $n\geq4$, by
\[
p_*(n)=1+\frac{\log(n-1)}{\log(n-2)}.
\]
Theorem~\ref{thm:star-equality} characterizes all extremizers, including the critical exponent and the distinct equality conditions at $p=\infty$. The equality proof for $p<1$ reduces a concave minimization problem to extreme points of finitely many polytopes and applies \cite{GKM}*{Proposition 6}. For $p>1$, we show that a maximizing function may be chosen with only one nonzero edge difference of the maximal function, reducing the sharp constant to a one-variable maximum.

Section~\ref{sec:prelim} contains the notation and auxiliary results, including the known complete-graph constant. Section~\ref{sec:universal} treats general finite graphs, and Section~\ref{sec:star} proves the sharp formula and the complete equality classification for star graphs.

\section{Preliminaries}\label{sec:prelim}
Throughout this paper, $G=(V,E)$ is a finite connected simple undirected graph, with $n=|V|\geq2$ and diameter $D$. We write $d_G$ for the graph distance and
\[
B_G(v,r)=\{w\in V:d_G(v,w)\leq r\}.
\]
All averages are taken with respect to counting measure. The centered Hardy--Littlewood maximal operator is defined by
\[
M_Gf(v)=\max_{0\leq r\leq D}\frac1{|B_G(v,r)|}
\sum_{w\in B_G(v,r)}|f(w)|,\qquad f:V\to\R,
\]
where the radii are integers. For $0<p<\infty$, define
\[
\Var_p(f)=\biggl(\sum_{\{u,v\}\in E}|f(u)-f(v)|^p\biggr)^{1/p},
\qquad
\Var_\infty(f)=\max_{\{u,v\}\in E}|f(u)-f(v)|,
\]
and put
\[
\C_{G,p}=\sup_{\Var_p(f)>0}\frac{\Var_p(M_Gf)}{\Var_p(f)},
\qquad 0<p\leq\infty.
\]
Each undirected edge is counted once. Since $G$ is connected, $\Var_p(f)=0$ if and only if $f$ is constant. We use $1/\infty=0$ and $t_+=\max\{t,0\}$.

The complete graph, star graph, and path graph on $n$ vertices are denoted by $K_n$, $S_n$, and $P_n$, respectively. We label the star vertices by $0,1,\ldots,n-1$, with center $0$. For $A\subset V$, let $\1_A$ be its indicator and let
\[
\partial_E A=\{\{u,v\}\in E:u\in A,\ v\notin A\}.
\]
We also write $\delta_w=\1_{\{w\}}$ and
$\osc(f)=\max_{v\in V}f(v)-\min_{v\in V}f(v)$.
For a finite vector $\mathbf z=(z_1,\ldots,z_N)$, we use the usual notation
\[
\|\mathbf z\|_p=\biggl(\sum_{j=1}^N|z_j|^p\biggr)^{1/p}\quad(0<p<\infty),
\qquad \|\mathbf z\|_\infty=\max_j|z_j|.
\]

The following elementary properties will be used without further comment.
\begin{lem}\label{lem:reduce}
The supremum defining $\C_{G,p}$ may be restricted to nonnegative functions. Moreover, for $t>0$,
\[
M_G(tf)=tM_Gf,\qquad \Var_p(tf)=t\Var_p(f),
\]
and $M_G(a+f)=a+M_Gf$ whenever $a\geq0$ and $f\geq0$.
\end{lem}
Indeed, $M_Gf=M_G|f|$ and $\Var_p(|f|)\leq\Var_p(f)$; the remaining identities follow from the ball averages.

\begin{lem}\label{lem:powers}
For $a_1,\ldots,a_k\geq0$ and $0<p\leq1$,
\[
\biggl(\sum_{i=1}^ka_i\biggr)^p
\leq\sum_{i=1}^ka_i^p
\leq k^{1-p}\biggl(\sum_{i=1}^ka_i\biggr)^p.
\]
Both inequalities reverse when $1\leq p<\infty$.
\end{lem}

\begin{lem}\label{lem:oscinput}
For $0<p\leq\infty$ and $f:V\to\R$,
\[
\osc(f)\leq D^{(1-1/p)_+}\Var_p(f).
\]
\end{lem}
This follows by summing edge differences along a shortest path from a minimum point to a maximum point and applying the finite-sequence inequality
$\|\mathbf z\|_1\leq D^{(1-1/p)_+}\|\mathbf z\|_p$.

\begin{lem}\label{lem:ordered-gaps}
Order the vertices as $v_1,\ldots,v_n$ so that
\[
t_j=f(v_j),\qquad t_1\leq\cdots\leq t_n,\qquad
\delta_j=t_{j+1}-t_j\quad(1\leq j<n).
\]
Then
\[
\sum_{j=1}^{n-1}\delta_j^p\leq\Var_p(f)^p\quad(1\leq p<\infty),
\qquad \max_j\delta_j\leq\Var_\infty(f).
\]
For $0<p\leq1$, one has $\sum_j\delta_j=\osc(f)\leq\Var_p(f)$.
\end{lem}
\begin{proof}
Let $A_j=\{v_1,\ldots,v_j\}$. Connectivity gives $|\partial_E A_j|\geq1$. For $1\leq p<\infty$,
\begin{align*}
\Var_p(f)^p
&=\sum_{\substack{a<b\\\{v_a,v_b\}\in E}}
\biggl(\sum_{j=a}^{b-1}\delta_j\biggr)^p\\
&\geq\sum_{\substack{a<b\\\{v_a,v_b\}\in E}}
\sum_{j=a}^{b-1}\delta_j^p
=\sum_{j=1}^{n-1}|\partial_E A_j|\delta_j^p
\geq\sum_{j=1}^{n-1}\delta_j^p.
\end{align*}
Every edge crossing $A_j$ has difference at least $\delta_j$, which proves the assertion for $p=\infty$. The assertion for $p\leq1$ follows from Lemma~\ref{lem:oscinput}.
\end{proof}

We will use the following standard extreme-point principle.
\begin{lem}\label{lem:vertex}
A continuous concave function on a nonempty compact convex polytope attains its minimum at an extreme point. More generally, a function that is continuous on a finite union of such polytopes and concave on each of them attains its minimum at an extreme point of one of the polytopes.
\end{lem}

The next inequality is due to Gonz\'alez-Riquelme, Kova\v{c}, and Madrid.
\begin{lem}[{\cite{GKM}*{Proposition 6}}]\label{lem:key}
Let $n\geq3$, $0<p<1$, and $2\leq r<n$ be an integer. Suppose that $x_1,\ldots,x_{n-r}\geq0$ and $y_1,\ldots,y_{r-1}\geq\eta\geq0$ satisfy
\[
\sum_{j=1}^{r-1}y_j=\sum_{i=1}^{n-r}x_i+n\eta.
\]
Then
\[
\sum_{i=1}^{n-r}x_i^p+(r-1)\eta^p
\leq\bigl(1-\tfrac1n\bigr)^p
\biggl(\sum_{i=1}^{n-r}x_i^p+\sum_{j=1}^{r-1}y_j^p\biggr).
\]
\end{lem}

We recall the sharp complete-graph constant for comparison. The finite-exponent result below is \cite{GKM}*{Theorem 1}; earlier ranges were obtained in \cite{GRM1}.
\begin{lem}\label{lem:complete-known}
For $n\geq2$ and $0<p\leq\infty$,
\[
\C_{K_n,p}=\frac{n-1}{n}.
\]
The indicator of any vertex attains equality.
\end{lem}
For $n=2$, the statement follows directly from the two ball averages. For $p=\infty$, use $\Var_\infty(f)=\osc(f)$ on $K_n$ and, for $f\geq0$,
\[
\osc(M_{K_n}f)=\max f-\frac1n\sum_{v\in V}f(v)
\leq\frac{n-1}{n}\osc(f).
\]
The indicator of one vertex gives the stated ratio for every exponent.

\section{Two-sided estimates for individual finite graphs}\label{sec:universal}
We first establish an optimal oscillation estimate and then combine it with point-mass tests and bounds for edge energy to obtain two-sided variational estimates.
Throughout this section, $n=|V|$, $D=\diam(G)$, and
\[
\gamma_p=(1-1/p)_+\quad(0<p\leq\infty),\qquad\gamma_\infty=1.
\]

\subsection{Oscillation estimates}\label{sec:osc}
\begin{lem}[Oscillation of the maximal function]\label{thm:osc}
For every nonnegative function $f$,
\[
\osc(M_Gf)\leq\bigl(1-\frac1n\bigr)\osc(f).
\]
\end{lem}
\begin{proof}
Every ball average is at most $\max f$, and a radius-zero ball at a maximum point attains this value. Every vertex can also use the full-graph ball. Hence
\[
\max_{v\in V}M_Gf(v)=\max_{v\in V}f(v),\qquad
\min_{v\in V}M_Gf(v)\geq\frac1n\sum_{v\in V}f(v).
\]
Choose $v_*$ with $f(v_*)=\max f$. Then
\begin{align*}
\osc(M_Gf)
&\leq\max_{v\in V}f(v)-\frac1n\sum_{v\in V}f(v)\\
&=\frac1n\sum_{v\in V\setminus\{v_*\}}
\bigl(\max_{u\in V}f(u)-f(v)\bigr)\\
&\leq\frac{n-1}{n}\osc(f).
\end{align*}
\end{proof}

\begin{thm}[The optimal uniform coefficient from variation to oscillation]\label{thm:sharp-osc-variation}
For $n\geq2$, let $\mathscr G_n$ be the class of connected simple undirected graphs on $n$ vertices. For $1<p<\infty$, let $p'=p/(p-1)$, so that $1/p+1/p'=1$. Define
\[
\Omega_p(n)=
\begin{cases}
1-1/n,&0<p\leq1,\\[2pt]
\displaystyle\frac1n\biggl(\sum_{j=1}^{n-1}j^{p'}\biggr)^{1/p'},&1<p<\infty,\\[4pt]
(n-1)/2,&p=\infty.
\end{cases}
\]
For every $G\in\mathscr G_n$ and every $f:V(G)\to\R$,
\[
\osc(M_Gf)\leq\Omega_p(n)\Var_p(f).
\]
Moreover, $\Omega_p(n)$ is the smallest coefficient valid for the whole class:
\[
\sup_{G\in\mathscr G_n}\ \sup_{\Var_p(f)>0}
\frac{\osc(M_Gf)}{\Var_p(f)}=\Omega_p(n).
\]
For every exponent, equality is attained on $P_n$.
\end{thm}

\begin{proof}
By Lemma~\ref{lem:reduce}, it suffices to prove the upper bound for $f\geq0$. Order its values as in Lemma~\ref{lem:ordered-gaps}. Every vertex can use the full-graph average, whereas $\max M_Gf=t_n$. Hence
\begin{equation}\label{eq:osc-weighted-gaps}
\osc(M_Gf)\leq t_n-\frac1n\sum_{i=1}^nt_i
=\frac1n\sum_{i=1}^{n-1}\sum_{j=i}^{n-1}\delta_j
=\frac1n\sum_{j=1}^{n-1}j\delta_j.
\end{equation}
For $1<p<\infty$, H\"older's inequality and Lemma~\ref{lem:ordered-gaps} give
\[
\sum_{j=1}^{n-1}j\delta_j
\leq\biggl(\sum_{j=1}^{n-1}j^{p'}\biggr)^{1/p'}
\biggl(\sum_{j=1}^{n-1}\delta_j^p\biggr)^{1/p}
\leq\biggl(\sum_{j=1}^{n-1}j^{p'}\biggr)^{1/p'}\Var_p(f).
\]
For the remaining exponents, the same lemma gives
\[
\sum_{j=1}^{n-1}j\delta_j\leq
\begin{cases}
(n-1)\Var_p(f),&0<p\leq1,\\
\dfrac{n(n-1)}2\Var_\infty(f),&p=\infty.
\end{cases}
\]
Together with \eqref{eq:osc-weighted-gaps}, these estimates prove the upper bound.

To prove sharpness, label the vertices of $P_n$ in path order as $v_1,\ldots,v_n$. Choose nonnegative gaps $\delta_1,\ldots,\delta_{n-1}$ and define
\[
f(v_1)=0,\qquad f(v_i)=\sum_{j=1}^{i-1}\delta_j\quad(2\leq i\leq n).
\]
The averages over the balls centered at $v_1$ are the averages of initial segments of this nondecreasing sequence. They are nondecreasing with the radius, so
\[
M_{P_n}f(v_1)=\frac1n\sum_{i=1}^nf(v_i),
\qquad M_{P_n}f(v_n)=f(v_n).
\]
These are respectively the minimum and maximum of $M_{P_n}f$. Consequently,
\[
\osc(M_{P_n}f)=\frac1n\sum_{j=1}^{n-1}j\delta_j,
\qquad
\Var_p(f)=\|(\delta_1,\ldots,\delta_{n-1})\|_p.
\]
For $1<p<\infty$, choose $\delta_j=j^{1/(p-1)}$. The ratio is then
\[
\frac{\osc(M_{P_n}f)}{\Var_p(f)}
=\frac{\sum_{j=1}^{n-1}j^{p'}}{n\bigl(\sum_{j=1}^{n-1}j^{p'}\bigr)^{1/p}}
=\Omega_p(n).
\]
For $p=\infty$, choose all gaps equal to one. For $0<p\leq1$, choose $\delta_{n-1}=1$ and all other gaps zero. Both choices give equality.
\end{proof}

\begin{rem}\label{rem:refined-osc-factor}
The preceding theorem and the diameter estimate can be combined by setting
\[
R_p(n,D)=\min\biggl\{D^{\gamma_p},\frac{\Omega_p(n)}{1-1/n}\biggr\}.
\]
Indeed, Lemmas~\ref{thm:osc} and \ref{lem:oscinput}, applied to $|f|$, give the diameter bound since $\Var_p(|f|)\leq\Var_p(f)$. Theorem~\ref{thm:sharp-osc-variation} gives the other bound. Therefore
\[
\osc(M_Gf)\leq\bigl(1-\frac1n\bigr) R_p(n,D)\Var_p(f).
\]
For $p\leq1$, $R_p(n,D)=1$, whereas $R_\infty(n,D)=\min\{D,n/2\}$.

\end{rem}

\subsection{Point-mass tests, diameter, and edge energy}
\begin{lem}[Improved lower bounds from point masses]\label{lem:simple-dirac-lower}
Suppose $n\geq4$. For every $0<p\leq\infty$, there is a vertex $w$ such that
\[
\frac{\Var_p(M_G\delta_w)}{\Var_p(\delta_w)}\geq\frac23.
\]
For $0<p\leq1$, the vertex can be chosen so that
\begin{equation}\label{eq:point-lower}
\frac{\Var_p(M_G\delta_w)}{\Var_p(\delta_w)}\geq\frac34.
\end{equation}
Here $\delta_w=\1_{\{w\}}$.
\end{lem}
\begin{proof}
For $v\ne w$, a ball centered at $v$ contains $w$ exactly when its radius is at least $d_G(v,w)$. Balls missing $w$ contribute zero. Once a ball contains $w$, its average is the reciprocal of its size, which cannot increase with the radius. Hence
\[
(M_G\delta_w)(w)=1,\qquad
(M_G\delta_w)(v)=|B_G(v,d_G(v,w))|^{-1}\quad(v\ne w).
\]
In particular, if $v$ is adjacent to $w$, then $(M_G\delta_w)(v)=(\deg(v)+1)^{-1}$. The only nonzero edge differences of the original function are the $\deg(w)$ edges incident to $w$, all with absolute value one. Thus
\[
\Var_p(\delta_w)=\deg(w)^{1/p}\quad(p<\infty),\qquad
\Var_\infty(\delta_w)=1.
\]
Connectivity guarantees $\deg(w)\geq1$.

\paragraph{The bound for all exponents.}
Write $\delta(G)=\min_{v\in V}\deg(v)$, and choose $w$ with $\deg(w)=\delta(G)$. If $\delta(G)\geq2$, each neighbor $v$ has degree at least $\delta(G)$, and
\[
|(M_G\delta_w)(w)-(M_G\delta_w)(v)|
=1-\frac1{\deg(v)+1}\geq1-\frac1{\delta(G)+1}\geq\frac23.
\]
For finite $p$, summing over these $\delta(G)$ edges gives $\Var_p(M_G\delta_w)^p\geq \delta(G)(2/3)^p$. Divide by $\delta(G)$ and take $p$th roots. For $p=\infty$, a single such edge suffices. If $\delta(G)=1$, the unique neighbor $v$ must have degree at least two; otherwise the pair $\{w,v\}$ would be a connected component, contradicting connectivity and $n\geq4$. The same edge then gives a difference at least $2/3$, and the variation of the original function is one for every exponent.

\paragraph{The improved bound for small exponents.}
Let $0<p\leq1$. Suppose a simple path $\mathcal P$ joins $w$ to $z$ and $(M_G\delta_w)(z)\leq1/4$. Denote its edge set by $\mathcal E(\mathcal P)$. Telescoping and Lemma~\ref{lem:powers} give
\begin{align*}
(3/4)^p
&\leq|(M_G\delta_w)(w)-(M_G\delta_w)(z)|^p\\
&\leq\biggl(\sum_{\{u,v\}\in\mathcal E(\mathcal P)}
|(M_G\delta_w)(u)-(M_G\delta_w)(v)|\biggr)^p\\
&\leq\sum_{\{u,v\}\in\mathcal E(\mathcal P)}
|(M_G\delta_w)(u)-(M_G\delta_w)(v)|^p.
\end{align*}
If $\mathcal P_1,\ldots,\mathcal P_{\delta(G)}$ are such paths and
$\mathcal E(\mathcal P_j)\cap\mathcal E(\mathcal P_k)=\varnothing$ for $j\ne k$, then
\[
\Var_p(M_G\delta_w)^p
\geq\sum_{j=1}^{\delta(G)}\ \sum_{\{u,v\}\in\mathcal E(\mathcal P_j)}
|(M_G\delta_w)(u)-(M_G\delta_w)(v)|^p
\geq \delta(G)(3/4)^p.
\]
We construct these paths by distinguishing the possible values of the minimum degree.

If $\delta(G)\geq3$, choose any minimum-degree vertex $w$. Each of its $\delta(G)$ neighbors satisfies $(M_G\delta_w)(v)\leq1/(\delta(G)+1)\leq1/4$, so the incident edges are the required paths.

If the minimum degree is one, take a leaf $w$ with neighbor $v$. If $\deg(v)\geq3$, use the edge $(w,v)$. Otherwise $\deg(v)=2$; let $u$ be its other neighbor. The vertex $u$ has a neighbor $z\ne v$, because otherwise $\{w,v,u\}$ would be a connected component. Since $w$ is a leaf, $z\ne w$. Thus the four vertices are distinct, $d_G(u,w)=2$, and
\[
\{w,v,u,z\}\subset B_G(u,2),\qquad
(M_G\delta_w)(u)\leq1/4.
\]
The path $(w,v,u)$ suffices.

Now suppose the minimum degree is two. If some vertex has degree at least three, follow a path from a degree-two vertex to such a vertex. Its first edge leaving the set of degree-two vertices can be denoted by $\{w,v\}$, with $\deg(w)=2$ and $\deg(v)\geq3$. Let $u$ be the other neighbor of $w$. The edge $(w,v)$ already ends at a vertex with maximal value at most $1/4$. If $\deg(u)\geq3$, use $(w,u)$ as the second path. If $\deg(u)=2$, let $z$ be its other neighbor. When $z=v$, the paths $(w,v)$ and $(w,u,v)$ have disjoint edge sets. When $z\ne v$, the vertex $z$ is not adjacent to $w$, whose only neighbors are $u,v$. Thus $d_G(z,w)=2$. Choose a neighbor $t\ne u$ of $z$, which exists by the minimum-degree assumption. Also $t\ne w$. Hence
\[
\{w,u,z,t\}\subset B_G(z,2),\qquad
(M_G\delta_w)(z)\leq1/4.
\]
The paths $(w,v)$ and $(w,u,z)$ have disjoint edge sets and both terminate at a vertex where $M_G\delta_w\leq1/4$.

Finally, if every vertex has degree two, $G$ is a cycle. Write the cycle order as $v_0,\ldots,v_{n-1}$. With $w=v_0$ and $z=v_2$, we have $d_G(w,z)=2$ and $B_G(z,2)$ contains at least $v_0,v_1,v_2,v_3$. Thus $(M_G\delta_w)(z)\leq1/4$. The two cycle arcs
\[
(v_0,v_1,v_2),\qquad
(v_0,v_{n-1},v_{n-2},\ldots,v_2)
\]
have disjoint edge sets and supply the two paths.
In every case, $w$ has minimum degree $\delta(G)$ and the sum of the $p$th powers of the maximal-function edge differences is at least $\delta(G)(3/4)^p$. Dividing by the corresponding sum $\delta(G)$ for the original function proves \eqref{eq:point-lower}.
\end{proof}

The next lemma combines the edge-count bound with a point-mass test at an endpoint of a diametral path.

\begin{lem}[Diameter, edge count, and a point-mass test]\label{lem:diam-simple}
Define
\[
N_E(n,D)=\min\biggl\{\frac{n(n-1)}2,
D+\frac{(n-D-1)(n-D+4)}2\biggr\}.
\]
Then $|E|\leq N_E(n,D)$. For every $0<p\leq\infty$, there is a vertex $w$ for which
\begin{equation}\label{eq:diameter-lower}
\frac{\Var_p(M_G\delta_w)}{\Var_p(\delta_w)}
\geq\frac{1-1/n}{(n-D)^{1/p}D^{\gamma_p}}.
\end{equation}
\end{lem}
\begin{proof}
The edge-count bound is due to Ore \cite{Ore}. For the second assertion, choose a shortest path $v_0,v_1,\ldots,v_D$ whose endpoints have distance $D$. Since $v_0$ is not adjacent to $v_2,\ldots,v_D$,
\[
\deg(v_0)\leq n-D,\qquad
\Var_p(\delta_{v_0})=\deg(v_0)^{1/p}\leq(n-D)^{1/p},
\]
with the degree factors interpreted as one at $p=\infty$. Moreover,
\[
(M_G\delta_{v_0})(v_0)=1,\qquad
(M_G\delta_{v_0})(v_D)=\frac1n.
\]
Indeed, the first ball centered at $v_D$ that contains $v_0$ has radius $D$ and equals $V$. Applying Lemma~\ref{lem:oscinput} to $M_G\delta_{v_0}$ gives
\[
1-\frac1n\leq\osc(M_G\delta_{v_0})
\leq D^{\gamma_p}\Var_p(M_G\delta_{v_0}).
\]
Combining the two variation estimates proves \eqref{eq:diameter-lower}.
\end{proof}

\begin{lem}[Edge energy controlled by oscillation]\label{lem:energy-size}
For $n\geq2$ and $0<p\leq\infty$, define
\[
Q_p(n,D)=
\begin{cases}
N_E(n,D)^{1-p}\min\{N_E(n,D),\lfloor n^2/4\rfloor\}^{p},&0<p<1,\\
\min\{N_E(n,D),\lfloor n^2/4\rfloor\},&1\leq p<\infty,\\
1,&p=\infty.
\end{cases}
\]
Every real-valued function $f$ satisfies
\[
\Var_p(f)\leq Q_p(n,D)^{1/p}\osc(f),
\]
where the factor at $p=\infty$ is one.
\end{lem}
\begin{proof}
The assertion is immediate for constant functions and for $p=\infty$. Otherwise, normalize $g=(f-\min f)/\osc(f)$, so $0\leq g\leq1$. For $1\leq p<\infty$, the convex function
\[
\Phi_p(\mathbf z)=\sum_{1\leq i<j\leq n}|z_i-z_j|^p
\]
attains its maximum on $[0,1]^n$ at a vertex. At a vertex with $k$ entries equal to one, its value is $k(n-k)\leq\lfloor n^2/4\rfloor$. Since $E$ is a subset of all vertex pairs,
\[
\Var_p(g)^p\leq\min\{|E|,\lfloor n^2/4\rfloor\}
\leq\min\{N_E(n,D),\lfloor n^2/4\rfloor\}=Q_p(n,D).
\]
For $0<p<1$, the finite-sequence inequality and the case $p=1$ give
\[
\Var_p(g)^p\leq |E|^{1-p}\Var_1(g)^p
\leq N_E(n,D)^{1-p}Q_1(n,D)^p=Q_p(n,D).
\]
Rescaling proves the result.
\end{proof}

\subsection{Two-sided variational estimates}\label{sec:general}
We now combine the preceding estimates. For $p=\infty$, we also use a maximum-degree test to improve the lower bound.

\begin{thm}[Two-sided bounds]\label{thm:generalbounds0}
Let $n\geq4$, and define
\[
b_p=\begin{cases}
3/4,&0<p\leq1,\\
2/3,&1<p<\infty,\\
3/4,&p=\infty,
\end{cases}
\]
\[
\ell_p(n,D)=\begin{cases}
\max\{b_p,\bigl(1-\frac1n\bigr)(n-D)^{-1/p}D^{-\gamma_p}\},&p<\infty,\\
\max\{3/4,1-n^{-1/D}\},&p=\infty.
\end{cases}
\]
Then
\begin{equation}\label{eq:general-bounds}
\ell_p(n,D)\leq\C_{G,p}
\leq\bigl(1-\frac1n\bigr) Q_p(n,D)^{1/p}R_p(n,D).
\end{equation}
\end{thm}
\begin{proof}
Lemma~\ref{lem:simple-dirac-lower} gives $\C_{G,p}\geq3/4$ for $0<p\leq1$ and $\C_{G,p}\geq2/3$ for all exponents. Lemma~\ref{lem:diam-simple} adds
\[
\C_{G,p}\geq\frac{\bigl(1-\frac1n\bigr)}{(n-D)^{1/p}D^{\gamma_p}}.
\]
This proves the stated finite-exponent lower bounds. We refine the infinite-exponent case next.

\paragraph{A uniform lower bound for $p=\infty$.}
Let $\Delta$ be the maximum degree. If $\Delta\geq3$, take a vertex $v$ of degree $\Delta$ and a neighbor $w$. Testing $\delta_w$ gives $\Var_\infty(\delta_w)=1$ and
\[
\Var_\infty(M_G\delta_w)
\geq|(M_G\delta_w)(w)-(M_G\delta_w)(v)|
=1-\frac1{\Delta+1}\geq\frac34.
\]
If $\Delta\leq2$, then $G$ is a path or a cycle.

On $P_n$, take $f(j)=j$ for $0\leq j\leq n-1$. Then $\Var_\infty(f)=1$. A ball is an integer interval, whose average equals the average of its endpoint values:
\[
\frac1{|B_{P_n}(j,r)|}\sum_{i\in B_{P_n}(j,r)}f(i)
=\frac{\max(0,j-r)+\min(n-1,j+r)}2.
\]
At $j=n-1$, this is at most $n-1$. At $j=n-2$ and $r\geq1$, the lower endpoint is at most $n-3$, so the average is at most $n-2$. Radius zero attains the central values in both cases. Consequently,
\[
M_{P_n}f(n-1)=n-1,\qquad M_{P_n}f(n-2)=n-2,
\]
and the variation of the maximal function is at least one.

For a cycle, write its vertices in cyclic order as $v_0,\ldots,v_{n-1}$ and set
\[
f(v_0)=0,\quad f(v_1)=1,\quad f(v_2)=2,\quad f(v_3)=1,
\qquad f(v_j)=0\quad(4\leq j<n).
\]
Then $\Var_\infty(f)=1$ and $\sum_{v\in V}f(v)=4$. At $v_2$, the singleton ball attains the maximum value two. At $v_1$ and $v_3$, the radius-one average is $(0+1+2)/3=1$, also when $n=4$. Balls of radius at least two have at least four vertices, so their averages are at most $4/4=1$. Radius zero also gives one. Thus
\[
M_Gf(v_2)=2,\qquad M_Gf(v_1)=M_Gf(v_3)=1.
\]
Again the variation of the maximal function is at least one. All cases therefore imply $\C_{G,\infty}\geq3/4$.

\paragraph{The diameter term for $p=\infty$.}
The same maximum-degree test, without assuming $\Delta\geq3$, gives
\[
\C_{G,\infty}\geq1-\frac1{\Delta+1}.
\]
Every vertex in $B_G(u,r+1)$ either belongs to $B_G(u,r)$ or is adjacent to one of its vertices. Hence
\[
B_G(u,r+1)\subset\bigcup_{z\in B_G(u,r)}
\bigl(\{z\}\cup\{y:\{z,y\}\in E\}\bigr),
\]
and
\[
|B_G(u,r+1)|\leq\sum_{z\in B_G(u,r)}(1+\deg(z))
\leq(\Delta+1)|B_G(u,r)|.
\]
Starting from $|B_G(u,0)|=1$ and iterating to radius $D$ gives $n\leq(\Delta+1)^D$. Thus $(\Delta+1)^{-1}\leq n^{-1/D}$, proving $\C_{G,\infty}\geq1-n^{-1/D}$.
This dominates the earlier term $\bigl(1-\frac1n\bigr)/D$. Indeed, with $\rho=n^{-1/D}\in(0,1)$,
\[
\frac{\bigl(1-\frac1n\bigr)}{D}=\frac{1-\rho^D}{D}
=\frac{1-\rho}{D}\sum_{j=0}^{D-1}\rho^j\leq1-\rho.
\]
Thus the definition of $\ell_\infty$ retains the stronger bound.

\paragraph{The upper bound.}
By Lemma~\ref{lem:reduce}, it suffices to consider $f\geq0$. Apply Lemma~\ref{lem:energy-size} to $M_Gf$ and then Remark~\ref{rem:refined-osc-factor}:
\begin{align*}
\Var_p(M_Gf)
&\leq Q_p(n,D)^{1/p}\osc(M_Gf)\\
&\leq\bigl(1-\frac1n\bigr) Q_p(n,D)^{1/p}R_p(n,D)\Var_p(f).
\end{align*}
Taking the supremum over $\Var_p(f)>0$ proves the upper bound. Each preceding test provides a lower bound for $\C_{G,p}$; their maximum is precisely the stated $\ell_p(n,D)$. The convention $1/\infty=0$ gives the same conclusion at $p=\infty$.
\end{proof}

The bounds in Theorem~\ref{thm:generalbounds0} need not be sharp for an individual graph.
In particular, the unboundedness result in \cite{GKM}*{Theorem 4} rules out an upper bound depending only on a fixed finite exponent $p$, uniformly in $n$.

\subsection{Bounds depending only on the number of vertices and the exponent}\label{sec:geometryfree}
Maximizing over all possible integer diameters removes the geometric parameter from the upper bound.

\begin{cor}[Bounds in terms of $n$ and $p$]\label{cor:general-geometryfree}
Let $n\geq4$, and put
\[
\mathcal H_p(n)=\max_{1\leq r<n}Q_p(n,r)^{1/p}R_p(n,r),
\]
where $r$ is an integer. Then
\[
b_p\leq\C_{G,p}\leq\bigl(1-\frac1n\bigr)\mathcal H_p(n).
\]
\end{cor}

\section{Sharp constants on star graphs}\label{sec:star}
On a star graph, the balls centered at a leaf are the singleton, the pair consisting of that leaf and the center, and the full graph. We therefore express the maximal function in terms of the signed center-to-leaf differences.
Gonz\'alez-Riquelme and Madrid \cite{GRM1} obtained the sharp constant for $p=1$, for all $p>1$ when $n=3$, and for several ranges with $0<p<1$. The constant and equality cases for $p=2$ and arbitrary $n\geq3$ were obtained in \cite{GRM2}. Theorem~\ref{thm:star} below gives the formula for every exponent; Theorem~\ref{thm:star-equality} describes every extremizer.

\subsection{Signed edge differences and the sharp formula}
Fix $S_n$, $n\geq3$, with center $0$ and leaves $1,\ldots,n-1$. By Lemma~\ref{lem:reduce}, upper bounds need only be proved for nonnegative functions. For such a function, write
\[
f(0)=x,\qquad f(i)=x+\xi_i\quad(1\leq i\leq n-1),
\qquad\boldsymbol\xi=(\xi_1,\ldots,\xi_{n-1}),
\]
and define
\begin{equation}\label{eq:star-s}
s=\frac1n\sum_{i=1}^{n-1}\xi_i.
\end{equation}
The coordinates $\xi_i=f(i)-f(0)$ are signed, and the full-graph average is
\[
\frac{f(0)+\sum_{i=1}^{n-1}f(i)}n=x+s.
\]
Thus $s=n^{-1}\sum_{\zeta\in V}f(\zeta)-f(0)$; this notation is used exclusively for the signed average in \eqref{eq:star-s}.

Conversely, every $\boldsymbol\xi\in\R^{n-1}$ is realized by a nonnegative function: choose
\[
x\geq\max\{0,-\xi_1,\ldots,-\xi_{n-1}\}.
\]
The translation identity in Lemma~\ref{lem:reduce} shows that the edge differences of $M_{S_n}f$ do not depend on the admissible choice of $x$.

\begin{lem}[Edge differences on a star]\label{lem:stard}
In the notation above, let $\eta_i=|M_{S_n}f(i)-M_{S_n}f(0)|$. Then
\[
\eta_i=
\begin{cases}
(\xi_i-s)_+,&s\geq0,\\
\xi_i,&s<0,\ \xi_i\geq0,\\
\min\{-\xi_i/2,-s\},&s<0,\ \xi_i<0.
\end{cases}
\]
With $\boldsymbol\eta=(\eta_1,\ldots,\eta_{n-1})$, one has
\[
\Var_p(f)=\|\boldsymbol\xi\|_p,\qquad
\Var_p(M_{S_n}f)=\|\boldsymbol\eta\|_p
\quad(0<p\leq\infty).
\]
\end{lem}
\begin{proof}
The balls centered at $0$ are $\{0\}$ and $V$; those centered at a leaf $i$ are $\{i\}$, $\{0,i\}$, and $V$. Thus
\[
M_{S_n}f(0)=x+\max\{0,s\},\qquad
M_{S_n}f(i)=x+\max\{\xi_i,\xi_i/2,s\}.
\]
If $s\geq0$, the difference is $(\xi_i-s)_+$. If $s<0$ and $\xi_i\geq0$, it is $\xi_i$. If $s<0$ and $\xi_i<0$, its absolute value is
\[
-\max\{\xi_i/2,s\}=\min\{-\xi_i/2,-s\}.
\]
The variation identities follow by summing over the edges $\{0,i\}$.
\end{proof}

The realization of every real edge vector and the reduction to nonnegative functions show that
\[
\C_{S_n,p}=\sup_{\boldsymbol\xi\ne\mathbf0}
\frac{\|\boldsymbol\eta(\boldsymbol\xi)\|_p}{\|\boldsymbol\xi\|_p}.
\]
We use this representation in both exponent ranges.

\begin{thm}[Sharp constants for every exponent]\label{thm:star}
For $n\geq3$, define
\[
p_*(3)=\infty,\qquad
p_*(n)=1+\frac{\log(n-1)}{\log(n-2)}\quad(n\geq4).
\]
For $1<p<\infty$, let $p'=p/(p-1)$. Then
\[
\C_{S_n,p}=
\begin{cases}
(n-1)/n,&0<p\leq1,\\[2pt]
\bigl((n-1)^{p'}+n-2\bigr)^{1/p'}/n,
&1<p<\infty,\ p\leq p_*(n),\\[2pt]
\bigl(1+(n-2)^{1-p}\bigr)^{-1/p},
&p_*(n)\leq p<\infty,\\[2pt]
1,&p=\infty.
\end{cases}
\]
For $n\geq4$, the two finite-exponent branches agree at $p=p_*(n)$. The infinite exponent is always covered by the last branch.
\end{thm}

\subsection{The known case \texorpdfstring{$0<p\leq1$}{0<p<=1}}
\begin{proof}[Proof of Theorem~\ref{thm:star} for $0<p\leq1$]
For $n\geq4$, the centered star-graph case of \cite{LLW}*{Theorem 1.3} gives $\C_{S_n,p}=1-1/n$. The case $n=3$ follows from \cite{LX}.
\end{proof}

\subsection{Proof for \texorpdfstring{$p>1$}{p>1}}
For $p>1$, we first obtain a lower bound from a family with one nonzero edge difference of the maximal function. We then show that every maximizer with $s>0$ belongs to this family up to a permutation and a positive scalar multiple. Vectors with $s\leq0$ are estimated separately.

\begin{proof}[Proof of Theorem~\ref{thm:star} for $p>1$]
Fix $1<p<\infty$, and retain the coordinates of Lemma~\ref{lem:stard}.

\paragraph{Step 1: a family giving the lower bound.}
For $0\leq t\leq1/(n-2)$, take
\[
\boldsymbol\xi_t=(1,-t,\ldots,-t),\qquad
s=\frac{1-(n-2)t}{n}\geq0.
\]
There are $n-2$ entries equal to $-t$. The first coordinate has edge difference of the maximal function $1-s$, and the other edge differences of the maximal function are zero. The resulting ratio is
\[
\mathcal R_p(t)=\frac{n-1+(n-2)t}{n\bigl(1+(n-2)t^p\bigr)^{1/p}}.
\]
It is realized by the nonnegative function
\[
f(0)=t,\qquad f(1)=1+t,\qquad f(i)=0\quad(2\leq i\leq n-1).
\]
Define
\[
\Gamma_p=\max_{0\leq t\leq1/(n-2)}\mathcal R_p(t).
\]
Continuity on this compact interval ensures that the maximum exists, and $\C_{S_n,p}\geq\Gamma_p$. The endpoint values give
\begin{equation}\label{eq:star-endpoint-bounds}
\Gamma_p\geq\bigl(1-\frac1n\bigr),\qquad
\Gamma_p\geq\bigl(1+(n-2)^{1-p}\bigr)^{-1/p}.
\end{equation}
We will prove that every edge vector has ratio at most $\Gamma_p$.

\paragraph{Step 2: vectors with $s\leq0$.}
If there is no positive coordinate, Lemma~\ref{lem:stard} gives $\eta_i\leq|\xi_i|/2$; when $s=0$, such a vector must be zero. Every nonzero vector of this type therefore has ratio at most $1/2<\bigl(1-\frac1n\bigr)\leq\Gamma_p$.

Suppose now that positive coordinates are present. List them as $u_i>0$ and the negative coordinates as $-\zeta_j<0$, $1\leq j\leq a$. Since the coordinate sum is nonpositive, at least one negative coordinate is needed, and the presence of a positive coordinate gives $1\leq a\leq n-2$. Set $S_+=\sum_i u_i>0$. Then
\[
\sum_{j=1}^a \zeta_j=S_+-ns.
\]
The edge differences of the maximal function at positive coordinates equal $u_i$, and each negative coordinate contributes at most $-s$. This remains true for $s=0$, when the maximal-function differences at negative coordinates vanish. Lemma~\ref{lem:powers} therefore yields
\begin{align*}
\sum_i \eta_i^p&\leq\sum_i u_i^p+a(-s)^p,\\
\sum_i|\xi_i|^p
&=\sum_i u_i^p+\sum_{j=1}^a \zeta_j^p
\geq\sum_i u_i^p+a^{1-p}(S_+-ns)^p.
\end{align*}
Since
\[
a^{1-p}(S_+-ns)^p\geq a^{1-p}n^p(-s)^p\geq a(-s)^p,
\]
the function
$U\mapsto[U+a(-s)^p]/[U+a^{1-p}(S_+-ns)^p]$ is nondecreasing for $U\geq0$.

Since $\sum_i u_i^p\leq S_+^p$ for $p>1$, we obtain
\begin{align*}
\frac{\sum_i \eta_i^p}{\sum_i|\xi_i|^p}
&\leq\frac{\sum_i u_i^p+a(-s)^p}
{\sum_i u_i^p+a^{1-p}(S_+-ns)^p}\\
&\leq\frac{S_+^p+a(-s)^p}
{S_+^p+a^{1-p}(S_+-ns)^p}\\
&\leq\frac{S_+^p+a(-s)^p}
{(1+a^{1-p})S_+^p+a^{1-p}n^p(-s)^p}.
\end{align*}
The last step uses $(S_++n(-s))^p\geq S_+^p+n^p(-s)^p$.

Moreover,
\begin{align*}
&(1+a^{1-p})S_+^p+a^{1-p}n^p(-s)^p\\
&\qquad=(1+a^{1-p})S_+^p+(n/a)^p\,a(-s)^p\\
&\qquad\geq\min\{1+a^{1-p},(n/a)^p\}
\bigl(S_+^p+a(-s)^p\bigr).
\end{align*}
Since $S_+>0$, the preceding quotient is bounded by
\[
\frac{\sum_i \eta_i^p}{\sum_i|\xi_i|^p}
\leq\max\biggl\{\frac1{1+a^{1-p}},\biggl(\frac an\biggr)^p\biggr\}
\leq\Gamma_p^p.
\]
For the last inequality, $1-p<0$ and $a\leq n-2$ imply
\[
\frac1{1+a^{1-p}}\leq\frac1{1+(n-2)^{1-p}}\leq\Gamma_p^p,
\qquad
(a/n)^p\leq\bigl(1-\frac1n\bigr)^p\leq\Gamma_p^p,
\]
by \eqref{eq:star-endpoint-bounds}. Taking $p$th roots completes Step 2.

\paragraph{Step 3: the shape of a maximizer with $s>0$.}
The ball-average formula in Lemma~\ref{lem:stard} expresses each $\eta_i$ as the absolute value of a difference of maxima of finitely many linear functions of $\boldsymbol\xi$. Thus $\boldsymbol\eta$ is continuous and satisfies $\boldsymbol\eta(t\boldsymbol\xi)=t\boldsymbol\eta(\boldsymbol\xi)$ for $t>0$. By homogeneity,
\[
\C_{S_n,p}=\max_{\boldsymbol\xi\in\Sigma_p}\|\boldsymbol\eta(\boldsymbol\xi)\|_p,
\qquad
\Sigma_p=\biggl\{\boldsymbol\xi\in\R^{n-1}:\sum_i|\xi_i|^p=1\biggr\}.
\]
The maximum is attained because $\Sigma_p$ is compact. It is positive by \eqref{eq:star-endpoint-bounds}. Fix a maximizing vector. If $s\leq0$, Step 2 gives $\C_{S_n,p}\leq\Gamma_p$; it remains to consider $s>0$.

On the open set of vectors satisfying $s>0$, define
\[
H(\boldsymbol\xi)=\sum_{i=1}^{n-1}(\xi_i-s)_+^p
-\C_{S_n,p}^p\sum_{i=1}^{n-1}|\xi_i|^p,
\qquad s=\frac1n\sum_i \xi_i.
\]
By the definition of $\C_{S_n,p}$, $H(\mathbf z)\leq0$ whenever $n^{-1}\sum_i z_i>0$, with equality at the chosen maximizing vector $\boldsymbol\xi$. Since this domain is open, $\boldsymbol\xi$ is an unconstrained local maximum of $H$.

The function $H$ is continuously differentiable for $p>1$. At the maximizing vector, let $\mathcal I=\{i:\xi_i>s\}$. This set is nonempty, since $\C_{S_n,p}>0$.

All partial derivatives of $H$ vanish at the local maximum. Dividing by $p$ gives, for active indices,
\begin{equation}\label{eq:star-stationary-active}
(\xi_j-s)^{p-1}-\frac1n\sum_{i\in\mathcal I}(\xi_i-s)^{p-1}
=\C_{S_n,p}^p \xi_j^{p-1}\qquad(j\in\mathcal I),
\end{equation}
and, for the other indices,
\begin{equation}\label{eq:star-stationary-inactive}
-\frac1n\sum_{i\in\mathcal I}(\xi_i-s)^{p-1}
=\C_{S_n,p}^p\sgn(\xi_j)|\xi_j|^{p-1}\qquad(j\notin\mathcal I).
\end{equation}
The left side of \eqref{eq:star-stationary-inactive} is strictly negative. Therefore every inactive coordinate is strictly negative, and
\[
|\xi_j|^{p-1}
=\frac{1}{n\C_{S_n,p}^p}
\sum_{i\in\mathcal I}(\xi_i-s)^{p-1}.
\]
The right-hand side is independent of $j$ and strictly positive. Thus all inactive coordinates are equal to the same negative number.

For $i\in\mathcal I$, we have $\xi_i>\xi_i-s>0$. Dividing \eqref{eq:star-stationary-active} by $\xi_i$ gives
\begin{align}
\C_{S_n,p}^p \xi_i^{p-2}
&=\frac{(\xi_i-s)^{p-1}-n^{-1}\sum_{j\in\mathcal I}(\xi_j-s)^{p-1}}{\xi_i}\notag\\
&<\frac{(\xi_i-s)^{p-1}}{\xi_i}
<(\xi_i-s)^{p-2}.\label{eq:star-second-compare}
\end{align}
The second strict inequality follows by writing its left side as
$(\xi_i-s)^{p-2}(\xi_i-s)/\xi_i$.
Since $\xi_i-s>0$, this comparison remains valid for $1<p<2$.

Suppose that $\mathcal I$ contains distinct indices $i,j$. Vary just these coordinates by
\[
\xi_i(\tau)=\xi_i+\tau,\qquad
\xi_j(\tau)=\xi_j-\tau,
\]
leaving the others fixed. Then
\[
\sum_{\ell=1}^{n-1}\xi_\ell(\tau)=\sum_{\ell=1}^{n-1}\xi_\ell,
\qquad
\xi_i(\tau)-s>0,\quad \xi_j(\tau)-s>0
\]
whenever $|\tau|<\frac12\min\{\xi_i-s,\xi_j-s\}$. All other coordinates and $s$ are fixed. The restriction of $H$ to this line is therefore twice differentiable near zero, and
\begin{align*}
\frac{d^2}{d\tau^2}H(\boldsymbol\xi(\tau))\bigg|_{\tau=0}
=p(p-1)\Bigl[&(\xi_i-s)^{p-2}+(\xi_j-s)^{p-2}\\
&-\C_{S_n,p}^p(\xi_i^{p-2}+\xi_j^{p-2})\Bigr]>0,
\end{align*}
by \eqref{eq:star-second-compare}. On the other hand, $H(\boldsymbol\xi(\tau))\leq H(\boldsymbol\xi)=0$ for all sufficiently small $|\tau|$, so its second derivative at zero is nonpositive. This contradiction proves $|\mathcal I|=1$.

After relabeling the leaves, the maximizing vector is $(a,-b,\ldots,-b)$, where $a,b>0$ and $a-(n-2)b=ns>0$. By homogeneity, we may divide this vector by $a$ and obtain
\[
(1,-t,\ldots,-t),\qquad t=b/a\in(0,1/(n-2)).
\]
Its ratio is $\mathcal R_p(t)\leq\Gamma_p$. Together with Step 2 and the lower bound in Step 1, this proves
\[
\C_{S_n,p}=\Gamma_p.
\]

\paragraph{Step 4: the one-variable maximum and extremizers.}
For $t>0$, logarithmic differentiation gives
\begin{align*}
\frac{\mathcal R_p'(t)}{\mathcal R_p(t)}
&=\frac{n-2}{n-1+(n-2)t}
-\frac{(n-2)t^{p-1}}{1+(n-2)t^p}\\
&=\frac{(n-2)\bigl[1-(n-1)t^{p-1}\bigr]}
{\bigl(n-1+(n-2)t\bigr)\bigl(1+(n-2)t^p\bigr)}.
\end{align*}
The denominator is positive. Put $t_0=(n-1)^{-1/(p-1)}$. The derivative is positive for $0<t<t_0$ and negative for $t>t_0$. By continuity at zero, the maximum on the allowed closed interval is therefore attained at
\[
t_p=\min\{(n-1)^{-1/(p-1)},1/(n-2)\}.
\]

If $t_0\leq1/(n-2)$, then $t_0^p=t_0/(n-1)=(n-1)^{-p'}$, and substitution gives
\[
\Gamma_p=\mathcal R_p(t_0)
=\frac{n-1}{n}\bigl(1+(n-2)(n-1)^{-p'}\bigr)^{1-1/p}
=\frac{\bigl((n-1)^{p'}+n-2\bigr)^{1/p'}}n.
\]

If $t_0\geq1/(n-2)$, the maximum is at the right endpoint, where
\[
\Gamma_p=\mathcal R_p(1/(n-2))
=\bigl(1+(n-2)^{1-p}\bigr)^{-1/p}.
\]
At $t_0=1/(n-2)$, the two expressions both equal $\mathcal R_p(1/(n-2))$ and therefore coincide.

For $n\geq4$, taking logarithms gives
\begin{align*}
t_0\leq\frac1{n-2}
&\quad\Longleftrightarrow\quad
(p-1)\log(n-2)\leq\log(n-1)\\
&\quad\Longleftrightarrow\quad p\leq p_*(n).
\end{align*}
For $n=3$, one has $t_0=2^{-1/(p-1)}<1$, so every finite $p>1$ falls in the first branch, consistently with $p_*(3)=\infty$.
The nonnegative extremizer in either branch is
\[
f(0)=t_p,\qquad f(1)=1+t_p,\qquad
f(i)=0\quad(2\leq i\leq n-1).
\]
Its edge vector is $\boldsymbol\xi_{t_p}$ and its ratio is $\mathcal R_p(t_p)=\Gamma_p$.

\paragraph{The infinite exponent.}
For $s\geq0$, $(\xi_i-s)_+\leq|\xi_i|$. For $s<0$, both other branches of Lemma~\ref{lem:stard} also give $\eta_i\leq|\xi_i|$. Hence
\[
\Var_\infty(M_{S_n}f)=\max_i \eta_i\leq\max_i|\xi_i|=\Var_\infty(f).
\]
For equality, take
\[
\boldsymbol\xi=(1,-1/(n-2),\ldots,-1/(n-2)).
\]
For this vector, $s=0$ and $\boldsymbol\eta=(1,0,\ldots,0)$. Since $n-2\geq1$, $\|\boldsymbol\xi\|_\infty=\|\boldsymbol\eta\|_\infty=1$. The function in Step 1 with $t=1/(n-2)$ realizes these differences, so $\C_{S_n,\infty}=1$.
\end{proof}

\subsection{All extremizers on star graphs}
We now characterize the nonconstant functions attaining $\C_{S_n,p}$ for every exponent.

\begin{thm}[Equality cases]\label{thm:star-equality}
Let $n\geq3$ and let $f:S_n\to\R$ be nonconstant.
\begin{enumerate}[label=\textup{(\roman*)}]
\item For $0<p<\infty$, put
\[
t_p=\begin{cases}
0,&0<p\leq1,\\[2pt]
\min\{(n-1)^{-1/(p-1)},1/(n-2)\},&1<p<\infty.
\end{cases}
\]
Then $f$ is an extremizer if and only if, after a permutation of the leaves, there are $a\geq0$, $u>0$, and $\varepsilon\in\{-1,1\}$ such that
\[
f(0)=\varepsilon(a+ut_p),\qquad
f(1)=\varepsilon\bigl(a+u(1+t_p)\bigr),\qquad
f(i)=\varepsilon a\quad(2\leq i\leq n-1).
\]
In particular, for $0<p\leq1$ these are precisely the positive scalar multiples of leaf indicators with a nonnegative constant added, together with their negatives.
\item For $p=\infty$, $f$ is an extremizer if and only if
\[
\frac1n\sum_{i=0}^{n-1}|f(i)|\leq|f(0)|,
\qquad
\max_{1\leq i\leq n-1}|f(i)|=|f(0)|+\Var_\infty(f).
\]
\end{enumerate}
\end{thm}
\begin{proof}
We first classify nonnegative extremizers for finite $p$. Retain $\xi_i$, $s$, $\eta_i$ from Lemma~\ref{lem:stard}. We then pass to arbitrary real functions and treat $p=\infty$ separately.

\paragraph{The range $0<p<1$ with $s\geq0$.}
Define
\[
\Psi(\boldsymbol\xi)=\bigl(1-\tfrac1n\bigr)^p
\sum_{i=1}^{n-1}|\xi_i|^p-\sum_{i=1}^{n-1}\eta_i^p.
\]
Theorem~\ref{thm:star} gives $\Psi(\boldsymbol\xi)\geq0$ for every edge vector, with equality at an extremizer.
At equality there is an index $i$ with $\xi_i>s$, since otherwise Lemma~\ref{lem:stard} gives zero maximal variation.
For $u>s$, set
\[
q_s(u)=\bigl(1-\tfrac1n\bigr)^pu^p-(u-s)^p.
\]
Since $p-1<0$ and $0<u-s\leq u$,
\[
q_s'(u)=p\biggl[\bigl(1-\tfrac1n\bigr)^pu^{p-1}-(u-s)^{p-1}\biggr]
\leq p\biggl[\bigl(1-\tfrac1n\bigr)^p-1\biggr]u^{p-1}<0.
\]
If $0<\xi_j\leq s$, replacing $(\xi_i,\xi_j)$ by
$(\xi_i+2\xi_j,-\xi_j)$ preserves their sum, and hence preserves $s$.
The $j$th coordinate contributes the same amount to $\sum|\xi_i|^p$ and zero to $\sum\eta_i^p$ before and after the replacement. Thus the change in $\Psi$ is
$q_s(\xi_i+2\xi_j)-q_s(\xi_i)<0$, contradicting its nonnegativity.
Consequently, every coordinate not exceeding $s$ is nonpositive.

After relabeling, write $\xi_i=a_i+s$ for $1\leq i\leq k$, where $a_i>0$ and $1\leq k\leq n-1$. The other coordinates are nonpositive, and
\[
\sum_{j=k+1}^{n-1}|\xi_j|=\sum_{i=1}^k a_i-(n-k)s\geq0.
\]
By Lemma~\ref{lem:powers},
\[
\sum_{i=1}^{n-1}|\xi_i|^p\geq F(s),\qquad
\sum_{i=1}^{n-1}\eta_i^p=\sum_{i=1}^k a_i^p,
\]
where
\[
F(t)=\sum_{i=1}^k(a_i+t)^p+\biggl(\sum_{i=1}^k a_i-(n-k)t\biggr)^p,
\quad 0\leq t\leq\frac{\sum_i a_i}{n-k}.
\]
Empty sums are zero, so this comparison also applies when $k=n-1$.
The function $F$ is continuous and strictly concave: in the interior,
\[
F''(t)=p(p-1)\biggl[\sum_{i=1}^k(a_i+t)^{p-2}
+(n-k)^2\biggl(\sum_i a_i-(n-k)t\biggr)^{p-2}\biggr]<0.
\]
At its left endpoint, Lemma~\ref{lem:powers} gives
$F(0)\geq(1+k^{p-1})\sum_i a_i^p\geq(1+1/k)\sum_i a_i^p$. Hence
\[
\bigl(1-\tfrac1n\bigr)^pF(0)>
\bigl(1-\tfrac1n\bigr)(1+1/k)\sum_i a_i^p\geq\sum_i a_i^p.
\]
At the right endpoint, for each $i$,
\[
\bigl(1-\tfrac1n\bigr)\biggl(a_i+\frac{\sum_j a_j}{n-k}\biggr)
\geq\bigl(1-\tfrac1n\bigr)\biggl(1+\frac1{n-k}\biggr)a_i\geq a_i.
\]
If $k\geq2$, the first inequality is strict, since $\sum_j a_j>a_i$.
If $k=1$, both inequalities are equalities. Comparing $F$ with the chord joining its endpoint values rules out every interior value of $s$ and the left endpoint. Equality therefore requires
\[
k=1,\qquad s=\frac{a_1}{n-1}.
\]
The total magnitude of the nonpositive coordinates is then
$a_1-(n-1)s=0$. Thus, after relabeling, $\boldsymbol\xi=(u,0,\ldots,0)$ with $u>0$.

\paragraph{The range $0<p<1$ with $s<0$.}
We show that $\Psi(\boldsymbol\xi)>0$ in this case.
If there are no positive coordinates, Lemma~\ref{lem:stard} gives
$\eta_i\leq|\xi_i|/2$, and the ratio is at most $1/2<1-1/n$.
Otherwise list the positive coordinates as $u_1,\ldots,u_k$ and the negative coordinates as $-\zeta_1,\ldots,-\zeta_a$, omitting zero coordinates. Set $S_+=\sum_i u_i>0$. Then
\[
k\geq1,\quad a\geq1,\quad k+a\leq n-1,\qquad
\sum_{j=1}^a\zeta_j=S_+-ns.
\]
By the edge formula,
\[
\Psi(\boldsymbol\xi)=\biggl[\bigl(1-\tfrac1n\bigr)^p-1\biggr]\sum_i u_i^p
+\sum_j g_s(\zeta_j),\qquad
g_s(\zeta)=\bigl(1-\tfrac1n\bigr)^p\zeta^p-\min\{\zeta/2,-s\}^p.
\]

If $k\leq(n-2)/2$, then
\[
g_s(\zeta)\geq\biggl[\bigl(1-\tfrac1n\bigr)^p-2^{-p}\biggr]\zeta^p,
\quad
\sum_j\zeta_j^p\geq(S_+-ns)^p>S_+^p\geq k^{p-1}\sum_i u_i^p.
\]
The first coefficient is positive. Moreover, since $pt^{p-1}$ is decreasing,
\begin{align*}
k^{1-p}\bigl(n^p-(n-1)^p\bigr)
&\leq kp(n-1)^{p-1}\\
&\leq\frac{n-2}{2}p(n-1)^{p-1}
\leq\int_{n/2}^{n-1}pt^{p-1}\,dt.
\end{align*}
Here $k^{1-p}\leq k$ because $k\geq1$. Dividing by $n^pk^{1-p}$ and rearranging gives
\[
\biggl[\bigl(1-\tfrac1n\bigr)^p-2^{-p}\biggr]k^{p-1}
\geq1-\bigl(1-\tfrac1n\bigr)^p.
\]
The preceding strict inequality therefore implies $\Psi(\boldsymbol\xi)>0$.

It remains to consider $k>(n-2)/2$. Now $a<n/2$, so
$S_+-ns>-2as$. Fix the positive coordinates and $s$, and minimize
$\sum_jg_s(\zeta_j)$ over the compact simplex
\[
\mathcal S_-=\biggl\{\boldsymbol\zeta\in[0,\infty)^a:
\sum_j\zeta_j=S_+-ns\biggr\}.
\]
For $J\subset\{1,\ldots,a\}$, let
\[
P_J=\biggl\{\boldsymbol\zeta\in\mathcal S_-:
0\leq\zeta_j\leq-2s\ (j\notin J),\quad\zeta_j\geq-2s\ (j\in J)\biggr\}.
\]
These sets cover the simplex. Each nonempty $P_J$ is a compact convex polytope, and the sum is concave on $P_J$, since
\[
g_s(\zeta)=
\begin{cases}
\biggl[\bigl(1-\tfrac1n\bigr)^p-2^{-p}\biggr]\zeta^p,&0\leq\zeta\leq-2s,\\
\bigl(1-\tfrac1n\bigr)^p\zeta^p-(-s)^p,&\zeta\geq-2s.
\end{cases}
\]
The minimum is therefore attained at an extreme point of one of these polytopes.
At such a point, at most one coordinate is outside $\{0,-2s\}$.
Indeed, two such coordinates could be changed by $+\delta,-\delta$ and by $-\delta,+\delta$, with $\delta>0$ smaller than their distances to $0$ and $-2s$. The two resulting vectors would remain in the same polytope and have the original vector as their midpoint, contradicting extremality.
Since the coordinate sum exceeds $-2as$, one coordinate is greater than $-2s$; all other nonzero coordinates equal $-2s$.

Write the nonzero coordinates of a minimizing vector as $y_1,\ldots,y_b$.
Then
\[
1\leq b\leq a,\qquad k+b\leq n-1,\qquad
y_j\geq-2s,\qquad\sum_{j=1}^b y_j=S_+-ns.
\]
Pad $u_1,\ldots,u_k$ with zeros to obtain $n-b-1$ nonnegative numbers.
All hypotheses of Lemma~\ref{lem:key} hold with $r=b+1$ and its parameter $\eta=-s$: in particular,
$2\leq b+1\leq n-k<n$ and $y_j\geq-s$.
The lemma gives
\[
\sum_i u_i^p+b(-s)^p
\leq\bigl(1-\tfrac1n\bigr)^p\biggl(\sum_i u_i^p+\sum_j y_j^p\biggr).
\]
We require strictness here. If $b\geq2$, change two coordinates to
$y_1\pm\delta,y_2\mp\delta$, with $0<\delta<-s$. Both perturbed vectors still have every coordinate at least $-s$ and have the same sum. Strict concavity of $\zeta^p$ implies that one perturbation strictly decreases $\sum_jy_j^p$. Applying Lemma~\ref{lem:key} to that perturbation preserves the left side and strictly decreases the right side. Thus the displayed inequality is strict for the original minimizing vector.

For $b=1$, the sole negative magnitude is $S_+-ns$. Since there are at most $n-2$ positive coordinates, concavity gives
\[
\sum_i u_i^p\leq(n-2)^{1-p}S_+^p.
\]
We use the strict numerical inequality
\begin{equation}\label{eq:strict-small-p}
(n-2)\bigl(n^p-(n-1)^p\bigr)^{1/(1-p)}+1
<(n-1)^{p/(1-p)}.
\end{equation}
Let $r=1/(1-p)>1$. The strictly decreasing derivative of $\zeta^p$ gives
\[
n^p-(n-1)^p=\int_{n-1}^n p\zeta^{p-1}\,d\zeta
<p(n-1)^{p-1}.
\]
Also $r-1=p/(1-p)>p^r$, and Bernoulli's inequality gives
\[
(n-1)^r\geq1+r(n-2)>n-1+p^r(n-2).
\]
Raising the first strict inequality to power $r$, multiplying by $n-2$, and combining with the last display proves \eqref{eq:strict-small-p}.
H\"older's inequality with exponents $1/p$ and $1/(1-p)$ now yields
\begin{align*}
&(\bigl(1-\frac1n\bigr)^{-p}-1)(n-2)^{1-p}S_+^p+\bigl(1-\frac1n\bigr)^{-p}(-s)^p\\
&\quad\leq(S_+-ns)^p
\biggl[(\bigl(1-\frac1n\bigr)^{-p}-1)^{1/(1-p)}(n-2)
 +(n-1)^{-p/(1-p)}\biggr]^{1-p}\\
&\quad<(S_+-ns)^p.
\end{align*}
Here $\bigl(1-\frac1n\bigr)^{-p}(-s)^p=(n-1)^{-p}(-ns)^p$, and the bracket is the left side of \eqref{eq:strict-small-p} divided by $(n-1)^{p/(1-p)}$.
Multiplying by $\bigl(1-\frac1n\bigr)^p$ and using the bound for $\sum_i u_i^p$, we obtain
\[
(1-\bigl(1-\frac1n\bigr)^p)\sum_i u_i^p+(-s)^p
<\bigl(1-\frac1n\bigr)^p(S_+-ns)^p.
\]
This is precisely the strict inequality needed when $b=1$.
Thus the minimum over the negative coordinates has strictly positive deficit also when $k>(n-2)/2$. No nonzero vector with $s<0$ is an extremizer for $0<p<1$.

\paragraph{The exponent $p=1$.}
Let $P=\sum_{\xi_i>0}\xi_i$ and $N=\sum_{\xi_i<0}|\xi_i|$, so $ns=P-N$. Suppose first that $s\geq0$ and the maximal function is nonconstant, and let $k=|\{i:\xi_i>s\}|\geq1$. Direct expansion gives
\[
\bigl(1-\frac1n\bigr)\sum_i|\xi_i|-\sum_i \eta_i
=\frac{(k-1)P+(n-1-k)N}{n}
 +\sum_{0<\xi_i\leq s}\xi_i.
\]
All terms are nonnegative. Since $P>0$, equality requires $k=1$, and then $n-1-k=n-2>0$ forces $N=0$. The last sum must also vanish. Hence exactly one coordinate is positive and all the others are zero.

For $s<0$, a vector without positive coordinates has ratio at most $1/2$. If $P>0$, let $a$ be the number of negative coordinates, so $1\leq a\leq n-2$. Since each negative coordinate contributes at most $-s=(N-P)/n$ to the maximal variation,
\[
\bigl(1-\frac1n\bigr)(P+N)-\sum_i \eta_i
\geq\frac{(a-1)P+(n-1-a)N}{n}>0.
\]
The strict sign follows from $N>P>0$ and $n-1-a\geq1$. This completes the classification for $0<p\leq1$.

\paragraph{The range $1<p<\infty$.}
For $s<0$ and at least one positive coordinate, Step 2 in the proof of Theorem~\ref{thm:star} is strict at the inequality
\[
(S_++n(-s))^p>S_+^p+n^p(-s)^p,
\]
because both summands are positive. All denominators there are positive, so the variational ratio is strictly less than $\C_{S_n,p}$. Vectors without positive coordinates again have ratio at most $1/2$.

If $s=0$ and $\boldsymbol\xi\ne\mathbf0$, there are both positive and negative coordinates. Write $P>0$ for their common total magnitude and $a$ for the number of negative coordinates. Then
\[
\frac{\sum_i \eta_i^p}{\sum_i|\xi_i|^p}
\leq\frac{1}{1+a^{1-p}}
\leq\frac{1}{1+(n-2)^{1-p}}.
\]
For the first inequality, the sum of the positive $p$th powers is at most $P^p$, whereas the sum of the negative $p$th powers is at least $a^{1-p}P^p$. Equality holds in both comparisons exactly when there is one positive coordinate, $a=n-2$, and all negative coordinates coincide. After scaling and relabeling, the vector is
\[
(1,-1/(n-2),\ldots,-1/(n-2)).
\]
The last bound equals $\C_{S_n,p}^p$ precisely when $n\geq4$ and $p\geq p_*(n)$, by the strict monotonicity on either side of the maximum of $\mathcal R_p$ established in Step 4. For $n=3$ and finite $p$, that maximum is attained strictly before the endpoint, so the endpoint bound is smaller.

If $s>0$, the stationary equations and the two-coordinate perturbation in Step 3 apply to every maximizing vector. They force, up to scaling and permutation,
\[
\boldsymbol\xi=(1,-t,\ldots,-t),\qquad 0<t<1/(n-2).
\]
The derivative of $\mathcal R_p$ has exactly one zero on $(0,\infty)$, at $(n-1)^{-1/(p-1)}$. Hence the only possible value is $t=t_p$, and it lies in the indicated open interval precisely when $p<p_*(n)$. At $p=p_*(n)$ the unique maximizing parameter is the endpoint and $s=0$; no additional finite-exponent equality cases appear there.

Consequently, every nonnegative finite-exponent extremizer has edge vector
$u(1,-t_p,\ldots,-t_p)$ with $u>0$. Nonnegativity of the other leaves gives $f(0)\geq ut_p$. Writing $a=f(0)-ut_p\geq0$ yields the values in \textup{(i)} with $\varepsilon=1$. Conversely, these values attain the sharp constant by the calculation of $\mathcal R_p(t_p)$, or by the leaf-indicator calculation when $p\leq1$. Adding $a\geq0$ and multiplying by $u>0$ preserve the variational ratio by Lemma~\ref{lem:reduce}.

\paragraph{Real-valued functions for finite $p$.}
For a nonconstant real function, equality in
\[
\Var_p(M_{S_n}f)=\Var_p(M_{S_n}|f|)
\leq\C_{S_n,p}\Var_p(|f|)
\leq\C_{S_n,p}\Var_p(f)
\]
requires both that $|f|$ be a nonnegative extremizer and that every edge satisfy
\[
|f(i)-f(0)|=\bigl||f(i)|-|f(0)|\bigr|.
\]
Indeed, the constant is positive, and a strict inequality on any edge would strictly increase the finite sum of $p$th powers. The displayed edge equality is equivalent to $f(i)f(0)\geq0$. If $f(0)\ne0$, all nonzero values therefore have the same sign. If $f(0)=0$, the nonnegative classification forces $p\leq1$ and a single nonzero leaf, whose sign is arbitrary. In either case a common sign $\varepsilon$ gives exactly \textup{(i)}. Such a common sign preserves both variations, proving the converse for real functions.

\paragraph{The exponent $p=\infty$.}
Apply Lemma~\ref{lem:stard} to $|f|$, with
\[
\xi_i=|f(i)|-|f(0)|,\qquad
s=\frac1n\sum_{i=1}^{n-1}\xi_i.
\]
For each edge, $|\xi_i|\leq|f(i)-f(0)|\leq\Var_\infty(f)$. If $s>0$, every positive maximal-function edge difference is $(\xi_i-s)_+<\Var_\infty(f)$, so equality is impossible. If $s\leq0$, a coordinate $\xi_i<0$ contributes at most $|\xi_i|/2$ (and contributes zero when $s=0$). A coordinate $\xi_i\geq0$ contributes exactly $\xi_i$. Since $\Var_\infty(f)>0$, equality therefore holds exactly when
\[
s\leq0,\qquad \max_i \xi_i=\Var_\infty(f).
\]
These are the two conditions in \textup{(ii)}. This argument proves necessity and sufficiency without imposing a common sign on $f$.
\end{proof}

\begingroup
\raggedright

\endgroup
\end{document}